\documentclass{article}
\usepackage{amssymb}
\usepackage{amsmath, amsthm}

\newtheorem{theorem}{Theorem}[section]
\newtheorem{corollary}{Corollary}[section]

\newtheorem{proposition}{Proposition}[section]
\newtheorem{remark}{Remark}

\title{Geometry in the large of the kernel of Lichnerowicz Laplacians and its applications}

\author{Vladimir Rovenski\footnote{Department of Mathematics, University of Haifa, 
Mount Carmel, Haifa, 31905, Israel, E-mail address: vrovenski@univ.haifa.ac.il},
\ Sergey Stepanov\footnote{Department of Mathematics, Russian Institute for Scientific
and Technical Information of the Russian Academy of Sciences,
20, Usievicha street, 125190 Moscow, Russia, E-mail address:  s.e.stepanov@mail.ru}
\ and \ Irina Tsyganok\footnote{Department of Data Analysis and Financial Technologies,
Finance University, 49-55, Leningradsky Prospect, 125468 Moscow, Russia,
E-mail address: i.i.tsyganok@mail.ru}
}

\begin{document}

\date{}
\maketitle

\begin{abstract}
There are very few general theorems on the kernel of the well-known Lichnerowicz Laplacian.
In the present article we consider the geometry of the kernel of
this operator restricted to covariant (not necessarily symmetric or skew-symmetric) tensors.
Our approach is based on the analytical method, due to Bochner, of proving vanishing theorems for the null space of  Laplace operator.
In particular, we pay special attention to the kernel of the Lichnerowicz Laplacian on Riemannian symmetric spaces of compact and noncompact types.
In conclusion, we give some applications to the theories of infinitesimal Einstein deformations and the stability of Einstein manifolds.
\end{abstract}



\section{\large Introduction}

Let $(M,g)$ be an $n$-dimensional $(n\ge 2)$ Riemannian manifold. In this case, the vector bundle $\otimes^{{p}}\,T^*M$ of covariant ${p}$-tensors $(1\le {p}<\infty)$ over $M$ carries the well-known {Lichnerowicz Laplacian} defined by the {Weitzenb\"{o}ck decomposition formula} (see \cite[p.~27]{1}; \cite[p.~54]{2};
\cite[p.~388--389]{3})
\begin{equation}\label{GrindEQ-1-1}
 \Delta_{L}=\bar{\Delta} + \Re_{p},
\end{equation}
where $\bar{\Delta}$ is the {rough Laplacian}
(see \cite[p.~54]{2} and \cite[p.~27]{1}) and $\Re_{{p}}:\,\otimes^{p}\,T^*M\to \otimes^{p}\,T^*M$ is the {Weitzenb\"{o}ck curvature operator} (see \cite[pp.~343--345]{4}) that depends linearly in known way on the Riemannian curvature tensor $Rm$ and the Ricci tensor ${\rm Ric}$ of the metric $g$.
In the present paper we will consider a one-parameter family of Lichnerowicz type Laplacians,
defined by the formula $\Delta_L = \bar\Delta_L + c\Re_{p}$ for any $c\in\mathbb{R}$.
For such Laplacians, we will retain the former original name ``Lichnerowicz Laplacian".
In addition, the subscript $L$ for $\Delta$ will denote that it is the
Lichnerowicz Laplacian in accordance with above definition.
On a closed (i.e., compact without boundary) Riemannian manifold $(M, g)$ we have an orthogonal (with respect to the global scalar product on $\otimes^{{p}}\,T^*M$) decomposition
\begin{equation}\label{GrindEQ-1-2}
 C^{\infty}( \otimes^{p}\,T^*M)={\rm Ker}\,\Delta_{L}\oplus {\rm Im}\,\Delta_{L},
\end{equation}
and the Lichnerowicz Laplacian has discrete eigenvalues with finite multiplicities,
which in few cases have been computed.
In \eqref{GrindEQ-1-2}, the first component ${\rm Ker}\,\Delta_{L}$ of \eqref{GrindEQ-1-2} is the kernel of the Lichnerowicz Laplacian $\Delta_{L}$. Its smooth sections are called $\Delta_{L}$-\textit{harmonic} (see \cite[p.~104]{5}). In the present article we study $\Delta_{L}$-harmonic sections of the bundle $\otimes^{{p}}\,T^*M$ using the analytical method, due ro S.\,Bochner (see, for example, the monographs
\cite[Chapter~9]{4} and \cite{5,6,12}), of proving vanishing theorems for the null space of a Laplace operator admitting a Weitzenb\"{o}ck decomposition formula \eqref{GrindEQ-1-1}.
In particular, we pay special attention to the kernel of the Lichnerowicz Laplacian on Riemannian symmetric spaces of compact and noncompact types.
In~conclusion we give applications of our results to the theories of \textit{infinitesimal Einstein deformations} and the \textit{stability of Einstein manifolds}
(see \cite{3} and \cite[Chapter~12]{2}).

\smallskip

\textbf{Acknowledgments}. The work
of 2nd and 3d authors
was supported by Russian Foundation for Basis Research (projects No. 16-01-00756).

\section{\large The kernel of the Lichnerowicz Laplacian ac\-ting on covariant $p$-tensors}

There are very few general theorems on the kernel ${\rm Ker}\,\Delta_{L}$ of the Lichnerowicz Laplacian
$\Delta_{L}: C^{\infty}(\otimes^{p}\,T^*M)\to C^{\infty}(\otimes^{p}\,T^*M)$ for the case $p\ge 2$.
In this section we fill this gap.

\smallskip

\textbf{2.1.}
Let ${M}$ be a differentiable manifold of dimension $n$ and $g$ be a Riemannian metric on ${M}$ with its Levi-Civita connection $\nabla$. To the Riemannian manifold $(M,\,g)$ one can associate a number of natural elliptic differential operators, which arise from the geometric structure of $(M,\,g)$. Usually these operators act in the space $C^{\infty}(E)$ of smooth sections of some vector bundle $E\to M$ of rank $p$ over $(M,\,g)$.
Moreover, let the vector bundle $E\to M$ has an inner product and its compatible connection $D$ induced by the metric $g$ and the Levi-Civita connection $\nabla$ of $(M, g)$, respectively.
 The connection $D$ induces a differential operator
$D:\,C^{\infty}(\otimes^{p}\,T^*M)\to C^{\infty}(\otimes^{p+1}\,T^*M)$
for any $p\in\mathbb{N}$, given by
\begin{eqnarray*}
 D T(Y,X_1,\ldots,X_p) &=& D_Y T(X_1,\ldots,X_p) \\
 &=& Y (T(X_1,\ldots,X_p)) -\sum\nolimits_j T(X_1,\ldots, D_Y X_j\ldots,X_p).
\end{eqnarray*}
Its formal adjoint $D^*:\,C^{\infty}(\otimes^{p+1}\,T^*M)\to C^{\infty}(\otimes^{p}\,T^*M)$
is given by
\[
 D^* T(X_1,\ldots,X_p) = -\sum\nolimits_j D_{e_i} T(e_i,X_1,\ldots,X_p),
\]
where $(e_j)$ is any local orthonormal frame of $TM$.
 Following the general theory (see, for example, the monographs \cite[p.~54]{2}; \cite[p.~105]{5} and \cite[pp.~308--333]{7}) we define the elliptic differential operator $A$ acting on $C^{\infty}(E)$ by means of the so-called \textit{Weitzenb\"{o}ck decomposition formula}
(see also \cite[p.~54]{2} and \cite{H2015})
\begin{equation}\label{GrindEQ-2-1}
 A=\bar{\Delta} + {c}\,\Re_{p},
\end{equation}
where the first component $\bar{\Delta}={D}^{*} D=-{\rm trace}(D^2)$ is the \textit{rough} (\textit{Bochner}) \textit{Laplacian} of $D$ and $D^{*}$ is the ${L}^{2} $-adjoint of $D$
(see \cite[p.~386]{3} and \cite[p.~52--53]{2}). The second term of the decomposition formula \eqref{GrindEQ-2-1}
contains a nonzero constant $c$ and the \textit{Weitzenb\"{o}ck curvature operator} $\Re$ that is a smooth symmetric endomorphism of $E$, which is in turn related to the curvature ${R}^{{D}}$ of a connection~$D$.

A section $\xi \in C^{\infty}(E)$ is called $A$-\textit{harmonic} if $A\xi =0$ (see \cite[p.~104]{5}).
We~define the vector space of $A$-harmonic $C^{\infty}$-sections of the vector bundle $E\to M$ by the condition
\[
 {\rm Ker}\, A=\{\xi \in C^{\infty}(E): A\,\xi =0 \}
\]
and the vector space of $A$-harmonic ${L}^{q}(E)$-sections of the vector bundle $E\to M$ by the condition
\[
 L^{q}( {\rm Ker}\, A)=\{\,\xi \in {\rm Ker}\, A:\|\xi\|\in L^{q} (M)\,\} .
\]
Furthermore, $A$-harmonic sections satisfy the (strong) unique continuation property. In local coordinates, the condition $A\,\xi =\bar{\Delta}\,\xi +c\,\Re_{p}(\xi)=0$ becomes a system of $p$ elliptic differential equations satisfying the structural assumptions of Aronszajn-Cordes (see Appendix of the monograph \cite{5}).
Consequently, the following proposition holds (see \cite[p.~104]{5}).

\begin{proposition}\label{T-2.1}
Let $E\to M$ be a rank ${r}$ vector bundle over the connected Riemannian manifold $(M, g)$ and let $A$ be an elliptic differential operator acting on $C^{\infty}$-sections of the vector bundle $E\to M$ and satisfying \eqref{GrindEQ-2-1}.
Let $\xi\in{\rm Ker}\,A$ be an $A$-harmonic section of the vector bundle $E\to M$.
If $\xi$ has a zero of infinite order at some point $x\in M$, then $\xi$ vanishes identically on $M$.
\end{proposition}

\textbf{2.2.} An illustration of the construction \eqref{GrindEQ-2-1} is given by the well-known Lichnerowicz Laplacian $\Delta_{L}: C^{\infty}(\otimes^{p}\,T^*M)\to C^{\infty}(\otimes^{p}\,T^*M)$, where $\otimes^{{p}}\,T^*M$ is the vector bundle of covariant $p$-tensors on a Riemannian manifold $(M,\,g)$. The general definition of $\Delta_{L}$ can be found in \cite[p.~344]{4}.
The Lichnerowicz Laplacian $\Delta_{L}$ is defined for an arbitrary covariant tensor $T\in C^{\infty}(\otimes^{p}\,T^*M)$ by the equality
\begin{equation}\label{GrindEQ-2-4}
 \Delta_{L}{T}=\bar{\Delta}\, T + {c}\,\Re_{p} (T).
\end{equation}
In this case, the Weitzenb\"{o}ck curvature operator $\Re_{p}:\,\otimes^{p}\,T^*M\to \otimes^{p}\,T^*M$ is defined by the following equalities (see also \cite[p.~54]{2}):
\begin{equation}\label{GrindEQ-2-5}
 (\Re_{p}(T))_{i_1,\ldots, i_p}
 =\sum\nolimits_{a} R_{i_a j}\, T^{\quad\, j}_{i_1\ldots\ \ldots\, i_{p}}
 -2\sum\nolimits_{a<b} R_{j\, i_a k\, i_b}
 T^{\quad j\ \ \ k}_{i_1\ldots\ \ldots\ \ldots\, i_p},
\end{equation}
where $T_{i_{1},\ldots,i_{p}}$, $R_{ij}$ and $R_{ijkl}$ are components of the tensor
$T\in C^{\infty}(\otimes^{p}\,T^*M)$, the Ricci tensor
and the Riemannian curvature tensor, respectively.
These components are defined by the following identities:
\[
 T_{i_{1}\ldots i_{p}}=T( e_{i_{1}},\ldots, e_{i_{p}}),\quad
 R_{ijkl} =g_{im}\,R_{jkl}^{m},\quad
 R_{kl} =R_{kil}^{i},
\]
where $R(e_{j}, e_{l} )\,e_{k}=R^{i}_{kjl} e_{i}$ and ${g}_{im} =g( e_{i}, e_{m})$ for an orthonormal frame $\{e_{1},\ldots, e_{n}\}$ of $T_{x} M$ at an arbitrary point $x\in M$ and for any $i, j, k,\ldots=1, 2,\ldots,n$. Moreover, the following identity holds:
\[
 g(\Re_{p}(T), T')=g(\Re_{p}(T'), T)
\]
for any $T, T'\in \otimes^{p}\,T^*M$, see \cite[p.~27]{1}.

\begin{remark}\label{R-1}\rm
 In \cite{4}, the Weitzenb\"{o}ck curvature operator \eqref{GrindEQ-2-5} is presented by
\begin{equation}\label{GrindEQ-2-5b}
 (\Re_{p}(T))(X_1,\ldots, X_p) = \sum\limits_{a,j}(R(e_j,X_a)T)(\underbrace{X_1,\ldots, e_j}_{a},\ldots, X_p).
\end{equation}
From the well known formula for curvature on the $(p,l)$-tensor bundle,
keeping in mind that $R(Y,Z)(T(X_1,\ldots, X_p))=0$ for a $p$-tensor $T$, we obtain
\[
 (R(Y,Z)T)(X_1,\ldots, X_p) =
 -\sum\nolimits_{a}T(X_1,\ldots, R(Y,Z)X_{a},\ldots, X_p).
\]
Thus, \eqref{GrindEQ-2-5b} can be rewritten in the form,
which obviously coincides with the ordinary Lichnerowicz Laplacian \eqref{GrindEQ-2-5},
\begin{eqnarray*}
 &&(\Re_{p}(T))(X_1,\ldots, X_p) = -2\sum\limits_{j,a; b<a} T(\underbrace{X_1,\ldots, R(e_j,X_a)X_b}_{b},\underbrace{\ldots, e_j}_{a-b},\ldots, X_p)\\
 &&\qquad  -\sum\limits_{j,a}T(\underbrace{X_1,\ldots, R(e_j,X_a)e_j}_{a},\ldots, X_p)\\
 &&\qquad = -2\sum\limits_{j,k,a; b<a} R(e_j,X_a,e_k,X_b)\cdot T(\underbrace{X_1,\ldots, e_k}_{b}, \underbrace{\ldots, e_j}_{a-b},\ldots, X_p)\\
 &&\qquad +\sum\limits_{j,a}{\rm Ric}(e_j,X_a)\cdot T(\underbrace{X_1,\ldots, e_j}_{a},\ldots, X_p).
\end{eqnarray*}
For special case $p=0$, the Lichnerowicz Laplacian $\Delta_{L}$ is the ordinary Laplacian
$\bar{\Delta}=-\,{\rm div}\,\circ\nabla$ acting on $C^{\infty}$-functions.
If $c=1$, then from the formula \eqref{GrindEQ-2-4} we obtain the classic definition of the Lichnerowicz Laplacian \eqref{GrindEQ-1-1}. In addition, if $\Delta_L$ acts on the bundle $\Lambda^p M$ of covariant skew-symmetric $p$-tensors
$(1\le p\le n-1)$ over $M$, then from \eqref{GrindEQ-2-4} we also obtain the Weitzenb\"{o}ck decomposition formula of the Hogde Laplacian $\Delta_H$
(see \cite[p.~347]{31}).
On the other hand, if $c=-1$ and $\Delta_L$ acts on the bundle $S^p M$ of covariant symmetric $p$-tensors $(1\le p<\infty)$ over $M$.
Then from \eqref{GrindEQ-2-4} we obtain the Weitzenb\"{o}ck decomposition formula for the Sampson Laplacian $\Delta_S$
(see \cite[p.~147]{S} and \cite[p.~55]{18}).
For $p=0$, the kernel of $\Delta_{L}$ consists of \textit{harmonic functions}.
For the case $p = 1$, we get $(\Re_{1}(T))(X)=T({\rm Ric}(X))$, thus
the Lichnerowicz Laplacian $\Delta_{L}$ has the form
$\Delta_{L}=\bar{\Delta}+{\rm Ric}$. In this case, the operator $\Delta_{L}$ is the \textit{Hodge-de Rham Laplacian} $\Delta_{{H}}$ acting on one-forms. Therefore, the kernel of $\Delta_{L}$ consists of \textit{harmonic one-forms} on $(M,g)$. Moreover, if $M$ is a closed manifold then we have the orthogonal decomposition \eqref{GrindEQ-1-2}, where the dimension of $L^{2}( {\rm Ker}\,\Delta_{L})$ equals to the first Betti number of $({M,g})$, according to de Rham's theorem. On the other hand, if $(M,g)$ is a Riemannian complete manifold with nonnegative Ricci curvature, then $L^{2}({\rm Ker}\,\Delta_{L})$ consists of parallel one-forms on $({M,g})$. Furthermore, if the Ricci curvature is positive at some point of $(M,g)$ or the holonomy of $(M,g)$ is irreducible then the vector space $L^{2}({\rm Ker}\,\Delta_{L})$ is trivial (for the proof, see \cite[p.~666]{8}).
 In particular, for $c=-1$ and $p=1$ we obtain from (4) that $\Delta_L=\overline\Delta -{\rm Ric}$.
 The kernel of this Laplacian consists of \textit{infinitesimal harmonic transformations} (see \cite{SM}).
Therefore, we will not consider these well-known cases, and in our article we assume that $p\ge 2$.
\end{remark}

 We can formulate the following corollary of Proposition~\ref{T-2.1}.

\begin{corollary}
Let $\otimes^{p}\,T^*M$ be the vector bundle of covariant $p$-tensors $(p\ge 2)$ on a connected Riemannian manifold $(M,g)$ and $\Delta_{L}$ the Lichnerowicz Laplacian acting on $C^{\infty}$-sections of $\otimes^{p}\,T^*M$.
If a $\Delta_{L}$-harmonic section $T$ of $\otimes^{p}\,T^*M$ has a zero of infinite order at some point $x\in M$, then $T$ vanishes identically on~$(M,g)$.
\end{corollary}

 By direct calculations, from \eqref{GrindEQ-2-4} we obtain the \textit{Bochner-Weitzenb\"{o}ck formula}
\begin{eqnarray}\label{GrindEQ-2-6}
\nonumber
 \frac{1}{2}\,\Delta_{B}\,(\|T\,\|^{2}) =-g(\bar{\Delta}\, T,  T) +\left\|\,\nabla\, T\,\right\|^{2} \\
 = - g( \Delta_{L}\, T, T) + \|\nabla\, T\,\|^{2} +{c}\,g(\Re_{p} (T), T),
\end{eqnarray}
where $c\ne0$ and $\Delta_{{B}} ={\rm div}\circ{\rm grad}$ is the \textit{Beltrami Laplacian} on $C^{\infty}$-functions.

Recall that the Riemannian curvature tensor of $(M, g)$ defines a symmetric algebraic operator $\bar{R}{:\Lambda}^{2}(T_{x} M)\to {\Lambda}^{2}(T_{x} M)$ on the vector space $\Lambda^{2}(T_{x} M)$ of all $2$-forms over~tangent space $T_{x} M$ at an arbitrary point $x\in M$ (see \cite[pp.~82--83]{4}). This $\bar{R}$ is called the \textit{curvature operator} of $(M, g)$.
The eigenvalues $\Lambda_{\alpha}$ of the curvature operator
$\bar{R}$ are real numbers at each point $x\in M$. Then we can select the orthonormal frame $\{\Xi_{\alpha}\}$ for $\Lambda^{2}({T}_{{x}}^{*} M)$ at each point $x\in M$, which consists of eigenvectors for $\bar{R}$, i.e., $\bar{R}(\Xi_{\alpha})={\it \Lambda}_{\alpha}\Xi_{\alpha}$.
In this case, the quadratic form $g(\Re_{p}(T), T)$ can be represented in the following form (see \cite[p.~346]{4}):
\begin{equation}\label{GrindEQ-2-7}
 {g}(\Re_{p}(T_{x}),\, T_{x})=\sum\nolimits_{\alpha}{\it \Lambda}_{\alpha}\|\,\Xi_{\alpha}(T_{x})\|^{2}
\end{equation}
at each point $x\in M$.

\begin{remark}\rm
There are many articles devoted to the relationship between the behavior of the curvature operator $\bar{R}$ of a Riemannian manifold $(M, g)$ and some global characterization of it, such as its homotopy type, topological type, etc. (see, e.g., \cite[pp.~351, 353, 390]{4} and \cite{9,10}).
In particular, Riemannian manifolds with nonnegative curvature operator are classified in \cite[Theorem~10.3.7]{6}.
\end{remark}

Recall that all the eigenvalues of $\bar{R}$ are all real numbers at each point $x\in M$.
Thus, we say that $\bar{R}$ is nonnegative (resp., strictly positive), or simply $\bar{R}\ge 0$ (resp., $\bar{R}>0$), if all the eigenvalues of $\bar{R}$ are nonnegative (resp., strictly positive).

Let $T$ be a $\Delta_{L}$-harmonic section of $\otimes^{p}\,T^*M$, then from \eqref{GrindEQ-2-6} we obtain
\begin{equation}\label{GrindEQ-2-8}
 \frac{1}{2}\,\Delta_{B}\,(\|T\,\|^{2}) = \|\,\nabla\,T\,\|^{2} +{c}\,g(\Re_{p} (T),\, T).
\end{equation}
Therefore, if $c>0$ and $g(\Re_{p}(T),\, T)\ge 0$ at any point of a connected
Riemannian manifold,
then $\Delta_{B}(\|T\|^{2})\ge 0$, and as a result of this, $\|T\|^{2}$ is a nonnegative subharmonic function. In this case, the following local theorem holds.

\begin{theorem}\label{T-2.2}
Let ${U}$ be a connected open domain of a Riemannian manifold $(M,g)$ with positive semi-definite curvature operator $\bar{R}$ at any point of ${U}$ and $\Delta_{L}: C^{\infty}(\otimes^{p}\,T^*M)\to C^{\infty}(\otimes^{p}\,T^*M)$ the Lichnerowicz Laplacian with $c>0$ acting on $C^{\infty}$-sections of the bundle of covariant ${p}$-tensor fields $\otimes^{p}\,T^*M$ over $(M,g)$ for $p\ge 2$. If $T\in {\rm Ker}\,\Delta_{L}$ at any point of ${U}$ and the scalar function $\|T\|^{\,2}$ has a local maximum at some point of ${U}$, then $\|T\|^{\,2}$ is constant and $T$ is invariant under parallel translation in ${U}$.
In addition, if $\bar{R}\ge k>0$ at some point $x\in U$ and $T\in C^\infty(\Lambda_p M)$ for all $p\in\{1,\ldots,n-1\}$ then $T\equiv 0$.
\end{theorem}

\noindent\textbf{Proof}.
From \eqref{GrindEQ-2-7} we conclude that the sign of the quadratic form $g(\Re_{p}(T), T)$ is opposite to the sign of the \textit{curvature operator} $\bar{R}$ of a Riemannian manifold $(M, g)$.
In particular, the formula \eqref{GrindEQ-2-7} immediately shows that the quadratic form $g(\Re_{p}(T_{x}),\,T_{x})$ is nonnegative (resp., positive) when the curvature operator is nonnegative (resp., positive) at each point $x\in M$.

Taking into account the above, we conclude from \eqref{GrindEQ-2-8} that if the curvature operator $\bar{R}$ of $(M,g)$ is positive semi-definite at any point of a connected open domain ${U}\subset {M}$, then $\|T\|^{2}$ is a subharmonic scalar function on $U$. Therefore, proceeding from \eqref{GrindEQ-2-8} and using the \textit{Hopf maximum principle} (see \cite[p.~26]{6} and \cite{11}), we can conclude that if the curvature operator of $(M,g)$ is positive semi-definite at any point of a connected open domain ${U}\subset {M}$, then $\|T\|^{2}$ is a constant $C\ge0$ and $\nabla\,T=0$ in ${U}$. If $C>0$, then $T$ is nowhere zero.

Now, at a point $x\in U$, where the curvature operator $\bar R$ satisfies the inequality $\bar R\ge k>0$,
we have
\begin{equation}\label{GrindEQ-2-7b}
 g(\Re_{p}(T_{x}),\,T_{x})=\sum\nolimits_{\alpha}\Lambda_{\alpha}\|\,\Xi_{\alpha}(T_{x})\|^{2}
 \ge k\sum\nolimits_{\alpha}\Lambda_{\alpha}\|\,\Xi_{\alpha}(T_{x})\|^{2}\ge0.
\end{equation}
In this case, the left hand side of \eqref{GrindEQ-2-7b} is zero, while the right hand side would be nonnegative.
This contradiction shows that $\Xi_{\alpha}(T_{x}) = 0$ for all $\alpha$. In particular, for $T\in C^\infty(\Lambda^p M)$
we have $T_x = 0$ (see \cite[p.~352]{24}). Then $C = 0$ and hence $T \equiv 0$.
\hfill$\Box$

\smallskip

Let $(M,g)$ be a closed Riemannian manifold. Note that \eqref{GrindEQ-2-8} is globally defined. Then there exists a point $x\in M$, at which the function $\|T\,\|^{2}$ attains the global maximum. At the same time, let $\|T\,\|^{2}$ satisfies the condition $\Delta_{B}(\|T\,\|^{2}) \ge 0$ everywhere in $(M,g)$. In this case, we can use the \textit{Bochner maximum principle}, which we deduce from the Hopf maximum principle.
Namely, it is well known that an arbitrary subharmonic function on a closed Riemannian manifold is constant (see \cite[Theorem~2.2]{6}).
Thus, the following statement automatically follows from our Theorem~\ref{T-2.2}.

\begin{corollary}\label{C-2.3}
Let $(M,g)$ be a connected and closed Riemannian manifold with positive semi-definite curvature operator $\bar{R}$ at each point of $(M,g)$ and $\Delta_{L}: C^{\infty}(\otimes^{p}\,T^*M)\to C^{\infty}(\otimes^{p}\,T^*M)$ the Lichnerowicz Laplacian with $c>0$ acting on $C^{\infty}$-sections of the bundle of covariant ${p}$-tensor fields $\otimes^{p}\,T^*M$  over $({M,g})$ for $p\ge 2$. If $T\in {\rm Ker}\,\Delta_{L}$ at any point of $(M,g)$, then $\|T\,\|^{\; 2}$ is a constant function and $T$ is invariant under parallel translation.
In addition, if $\bar{R}\ge k>0$ at some point $x\in M$ and $T\in C^\infty(\Lambda_p M)$ for all $p\in\{1,\ldots,n-1\}$ then $T\equiv 0$.
\end{corollary}

\begin{remark}\rm
The well-known Hodge-de Rham Laplacian $\Delta_{H}$ acting on $C^{\infty}$-secti\-ons of the bundle of ${p}$-forms $\Lambda^{p}(M)$ is the most famous example of the Lichnerowicz Laplacian $\Delta_{L}$. An arbitrary $\Delta_{H}$-harmonic section of $\Lambda^{p} M$ is called a \textit{harmonic $p$-form} $(1\le p\le n-1)$. Moreover, by Hodge theory, the $p$-th Betti number of a closed manifold $M$ is precisely $\beta_{p}(M)={\rm dim}\,{\rm Ker}\,\Delta_{H} $, the dimension of the space of harmonic $p$-forms on $(M,g)$. Thus we conclude from our Theorem~\ref{T-2.2} that a closed Riemannian manifold with positive curvature operator has vanishing the ${p}$-th \textit{Betti number} $\beta_{p}(M)$. The added benefit is that we can also conclude from our theorem that if the curvature operator is merely nonnegative, then the ${p}$-th Betti number satisfies the inequality
$\beta_{p}(M)\le \Big(\begin{array}{c} {n} \\ {p}\end{array}\Big)$
(see the {Meyer-Gallot Theorem} in \cite[p.~351]{4}).
\end{remark}

As analogues of Theorem~\ref{T-2.2} and Corollary~\ref{C-2.3}, we can prove the following theorem and corollary.

\begin{theorem}
Let $U$ be a connected open domain of a Riemannian manifold $(M, g)$ with negative semi-definite curvature operator $\bar R$ at any point of $U$ and $\Delta_L: C^\infty\otimes^p\,T^* M \to C^\infty\otimes^p\,T^* M$ be the Lichnerowicz Laplacian with $c < 0$ acting on $C^\infty$-sections of the bundle of covariant $p$-tensor fields $\otimes^p\,T^* M$ over $(M,g)$ for $p\ge2$.
If $T\in{\rm Ker}\,\Delta_L$ at any point of $U$ and the scalar function $\|T\|^2$ has a local maximum at some point of $U$, then $\|T\|^2$ is a constant function and $T$ is invariant under parallel translation in $U$.
In addition, if $\bar{R}\le k<0$ at some point $x\in M$ and $T\in C^\infty(\Lambda_p M)$ for all $p\in\{1,\ldots,n-1\}$ then $T\equiv 0$.
\end{theorem}


\begin{corollary}
Let $(M, g)$ be a connected and closed Riemannian manifold with negative semi-definite curvature operator $R$ at each point and $\Delta_L: C^\infty\otimes^p\,T^* M \to C^\infty\otimes^p\,T^* M$ be the Lichnerowicz Laplacian with $c < 0$ acting on $C^\infty$-sections of the bundle of covariant $p$-tensor fields $\otimes^p\,T^* M$ over $(M,g)$ for $p\ge 2$.
If $T\in{\rm Ker}\,\Delta_L$ at any point of $(M, g)$, then $\|T\|^2$ is a constant function and $T$ is invariant under parallel translation.
In addition, if $\bar{R}\le k<0$ at some point $x\in M$ and $T\in C^\infty(\Lambda_p M)$ for all $p\in\{1,\ldots,n-1\}$ then $T\equiv 0$.
\end{corollary}

 By direct calculation we find the following:
\[
 \frac{1}{2}\,\Delta_{B}\,(\|T\|^{2}) = \|T\|\cdot\Delta_{B}(\|T \|+\|\, d\,\| T \|\,\|^{2} .
\]
Then the equation \eqref{GrindEQ-2-8} can be rewritten in the form
\[
 \|T\,\|\,\Delta_{B}(\|T\,\|) = {c}\,g(\Re_{p} (T),\, T)+\|\nabla\, T\,\|^{2}
 -\| d\,\|T\,\|\,\|^{2} .
\]
Using the \textit{Kato inequality} (see \cite{13})
\[
 \|\nabla\, T\,\|^{2} \ge \|\, d\,\|T\,\|\,\|^{2},
\]
we can write the following inequality (with $c\ne0$):
\begin{equation}\label{GrindEQ-2-9}
 \| T\|\cdot \Delta_{B}(\|T\|) \ge {c}\,g(\Re_{p}(T), T).
\end{equation}
Therefore, if we suppose that $c>0$ and
$g(\Re_{p} (T), T)\ge 0$ at any point of a Riemannian manifold $(M, g)$, then we have $\Delta_{B}\,(\|T\,\|) \ge 0$ and as a result of this $\|T\|$ is a nonnegative subharmonic function on $(M,g)$.

On the other hand, Greene and Wu proved in \cite{14} the following proposition:
If $(M, g)$ is a complete and noncompact Riemannian manifold with nonnegative sectional curvature and $f$ is a nonnegative subharmonic function on $(M, g)$,
then $\int_{M} f^{q}\,{\rm d}\,{\rm vol}_g =+\infty$ for any $1\le q<+\infty$ unless $f\equiv 0$. Based on \eqref{GrindEQ-2-9} and using our Theorem~\ref{T-2.2} and the Greene-Wu theorem, we conclude that if $T\in {\rm Ker}\,\Delta_{L}$ at any point of $(M, g)$ and $\int_{M}\|T\,\|^{\,q}\,{\rm d}\,{\rm vol}_g <+\infty$ for some $1\le q<+\infty $, then $T\equiv 0$ for the case of a complete noncompact Riemannian manifold $(M, g)$ with nonnegative sectional curvature. Therefore, we can formulate the following.

\begin{theorem}\label{T-2.3}
Let $(M,g)$ be a complete noncompact Riemannian manifold with positive semi-definite curvature operator $\bar{R}$ and $\Delta_{L}:\,C^{\infty}(\otimes^{p}\,T^*M)\to C^{\infty}(\otimes^{p}\,T^*M)$ be
the Lichnerowicz Laplacian with $c>0$, acting on $C^{\infty}$-sections of the bundle $\otimes^{p}\,T^*M$ of covariant $p$-tensor over $(M,g)$ for $p\ge 2$. Then the vector space $L^{q}( {\rm Ker}\,\Delta_{L})$ is trivial for an arbitrary $1\le q<+\infty$.
\end{theorem}

As an analogue of Theorem~\ref{T-2.3} we can prove the following.

\begin{theorem}\label{T-2.4}
Let $(M, g)$ be a complete simply connected Riemannian manifold with negative semi-definite curvature operator
$\bar R$ and $\Delta_{L}:\,C^{\infty}(\otimes^{p}\,T^*M)\to C^{\infty}(\otimes^{p}\,T^*M)$ the Lichnerowicz Laplacian with $c<0$ acting on $C^{\infty}$-sections of the bundle of covariant ${p}$-tensor fields $\otimes^{p}\,T^*M$ over $(M,g)$ for $p\ge 2$. If~$T\in L^q({\rm Ker}\,\Delta_{L})$
 for some $q\in(0,+\infty)$, then $\|T\|$ is a constant function and $T$ is invariant under parallel translation.
 In particular, if
 ${\rm vol}(M,g)= +\infty$, then $T\equiv0$.
 In addition, if $\bar{R}\le k<0$ at some point $x\in M$ and $T\in C^\infty(\Lambda^p M)$ for all $p\in\{1,\ldots,n-1\}$ then $T\equiv 0$.
\end{theorem}

\noindent\textbf{Proof}.
Note that if $\bar R\le 0$, then $g(\Re_{p}(T_{x}),\,T_{x})\le0$ at an arbitrary point $x\in M$ because
$g(\Re_{p}(T_{x}),\,T_{x})=\sum\nolimits_{\alpha}\Lambda_{\alpha}\|\,\Xi_{\alpha}(T_{x})\|^{2}\le0$,
where $\Lambda_{\alpha}\le0$ for all $\alpha$. On the other hand, if we suppose
that $c<0$ and $g(\Re_{p}(T_{x}),\,T_{x})\le0$ at any point of $(M,g)$, then from \eqref{GrindEQ-2-7b}
we obtain the inequality $\Delta_B(\|T\|\ge0$. In this case, we conclude that $\|T\|$ is a nonnegative subharmonic function on $(M,g)$.
On~the other hand, in \cite[p.~288]{LS-84} was proved that every nonnegative subharmonic $L^q$-function
for $q\in(0,+\infty)$ on a complete simply connected Riemannian manifold $(M,g)$ of nonpositive
sectional curvature is constant.
Therefore, if $(M,g)$ is a complete simply connected Riemannian manifold
and $T\in L^q({\rm Ker}\,\Delta_L)$ for some $q\in(0,+\infty)$, then $\|T\|$ is a constant $C\ge 0$ and $\nabla T=0$.

In this case, the inequality $\int_M \|T\|^q\,d\,{\rm vol}_g<+\infty$
can be rewritten in the form  $C^q\int_M d\,{\rm vol}_g=C^q{\rm vol}(M,g)<+\infty$
Therefore, if ${\rm vol}(M,g)=+\infty$ then $C=0$ and, hence, $T\equiv0$.
If $C>0$, then $T$ is nowhere zero. Now, at a point $x\in M$, where the curvature operator $\bar R$ satisfies the inequality $\bar R\le k<0$, we have
\[
  g(\Re_{p}(T_{x}),\,T_{x})=\sum\nolimits_{\alpha}\Lambda_{\alpha}\|\,\Xi_{\alpha}(T_{x})\|^{2}
 \le k\sum\nolimits_{\alpha}\Lambda_{\alpha}\|\,\Xi_{\alpha}(T_{x})\|^{2}\le0.
\]
In this case, the left hand side of \eqref{GrindEQ-2-7b} is zero, while the right hand side would be nonpositive.
This contradiction shows that $\Xi_{\alpha}(T_{x}) = 0$ for all $\alpha$. In particular, for $T\in C^\infty(\Lambda^p M)$ we have $T_x = 0$ (see \cite[p.~352]{24}). Then $C = 0$ and hence $T \equiv 0$.
\hfill$\Box$

\smallskip

\textbf{2.3.} Recall that the \textit{Riemannian symmetric space} is a finite dimensional Riemannian manifold $(M,g)$, such that for every its point $x$ there is an involutive geodesic symmetry $s_{x} $, such that $x$ is an isolated fixed point of $s_{x}$. $(M,g)$ is said to be \textit{Riemannian locally symmetric} if its geodesic symmetries are in fact isometries. A Riemannian locally symmetric space is said to be a \textit{Riemannian globally symmetric space} if, in addition, its geodesic symmetries are defined on all $(M,g)$.
A Riemannian globally symmetric space is complete (see \cite[p.~244]{15}).
In addition, a complete and simply connected~Riemannian locally symmetric~space is a Riemannian globally symmetric space (see~\cite[p.~244]{15}). Riemannian globally symmetric spaces can be classified in terms of their isometry groups.
The classification distinguishes three basic types of Riemannian globally symmetric spaces: the spaces of so-called \textit{compact type}, the spaces of so-called \textit{noncompact type} and the spaces of \textit{Euclidean type} (see, for example,
\cite[p.~252]{15}). An addition, if $(M,g)$ is a Riemannian globally symmetric space of compact type then $(M,g)$ is a closed Riemannian manifold with non-negative sectional curvature and positive-definite Ricci tensor
(see \cite[p.~256]{15}). Moreover, its curvature operator $\bar{R}$ is nonnegative (see \cite{16}). Using Remark~\ref{R-1} and Corollary~\ref{C-2.3}, one can argue that the following statement holds.

\begin{proposition}
Let $(M,g)$ be an $n$-dimensional $({n}\ge 2)$ simply connected Riemannian globally symmetric space of compact type and $\Delta_{L}: C^{\infty}(\otimes^{p} T^*M)$ $\to C^{\infty}(\otimes^{p} T^*M)$ the Lichnerowicz Laplacian acting on $C^{\infty}$-sections of the bundle $\otimes^{p}\,T^*M$ of covariant $p$-tensors over $(M,g)$ for $p\ge 1$. Then all $\Delta_{L}$-harmonic sections of $T^*M$ vanish everywhere on $(M,g)$ and every $\Delta_{L}$-harmonic section of $\otimes^{p}\,T^*M$ for the case $p\ge 2$ is invariant under parallel translation.
\end{proposition}

\begin{remark}\rm
If $(M,g)$ is a simply connected Riemannian globally symmetric space of compact type then its Betti numbers satisfy the following conditions: $b_{1}(M)=b_{n-1}(M)=0$ and $b_{p}(M)\le \Big(\begin{array}{c} {n} \\ {p} \end{array}\Big)$ for $p=2,\ldots, n-2$.
This statement follows directly from Theorem~\ref{T-2.3}.
\end{remark}

Recall that a Riemannian globally symmetric space $(M, g)$ is complete.
Moreover, a Riemannian symmetric space of
noncompact type has the nonpositive sectional curvature and negative-definite Ricci tensor, see \cite[p.~256]{33}.
It is also known that a Riemannian symmetric space has the nonpositive curvature operator if and only if it has the nonpositive sectional curvature (see \cite{16}). After the above remarks, the assertion of the following
proposition becomes an obvious corollary of Theorem~\ref{T-2.4}.

\begin{proposition}
Let $(M,g)$ be an $n$-dimensional $(n\ge 2)$ simply connected Ri\-emannian globally symmetric space of noncompact type and $\Delta_L:\,C^{\infty}(\otimes^{p}\,T^*M)$ $\to C^{\infty}(\otimes^{p}\,T^*M)$ be the Lichnerowicz Laplacian with $c < 0$ acting on
$C^\infty$-sections of the bundle $\otimes^{p}\,T^*M$ of covariant $p$-tensors over $(M,g)$ for $p\ge 2$.
Then the vector space $L^q({\rm Ker}\,\Delta_L)$ for any $q\in(0,+\infty)$ consists of parallel tensor fields.
In particular, if
${\rm vol}(M,g)= +\infty$, then
$L^q({\rm Ker}\,\Delta_L)$ is trivial.
\end{proposition}

\section{\large The kernel of the Lichnerowicz Laplacian ac\-ting on symmetric bilinear forms}

Here, we study the kernel ${\rm Ker}\,\Delta_{L}$ of the Lichnerowicz Laplacian
$\Delta_{L}: C^{\infty}(S^{p} M)$ $\to C^{\infty}(S^{p} M)$ restricted to the space of $C^{\infty}$-sections $C^{\infty}(S^{p} M)$ of the vector bundle $S^{p} M$ of covariant symmetric $p$-tensors $(p\ge 2)$ and, in particular, of $C^{\infty}$-sections $C^{\infty}(S^{2} M)$ of the vector bundle $S^{2} M$ of symmetric bilinear forms $S^{2} M$ on a Riemannian manifold $(M,g)$. In this section, we denote by $\varphi$ an arbitrary smooth section of $C^{\infty}(S^{p} M)$.

\smallskip

\textbf{3.1}. First, we consider the subbundle $S_{0}^{p} M$ of the bundle $\otimes^{p}\,T^*M$ consisting of smooth traceless symmetric tensor fields. A~section $\varphi\in C^{\infty}(S_{0}^{p} M)$ is defined by the condition
\[
 {\rm trace}_{g}\,\varphi :=\sum\nolimits_{i=1}^{n}\,\varphi( e_{i}, e_{i}, X_{3},\ldots, X_{p}) =0
\]
for the orthonormal frame $\{e_{i}\}$ of $T_x M$ at an arbitrary point $x\in M$.
It is well known (see \cite{3}) that $g(\Re_2(\varphi),\varphi)\ge0$ for any $\varphi\in C^\infty(S^2_0 M)$ if ${\rm sec}\ge0$ for the sectional curvature of $(M,g)$.
This statement was generalized in \cite[p.~8]{17} in the following form: $g(\Re_p(\varphi),\varphi)$ is positive-semidefinite for $p\ge2$ if ${\rm sec}\ge0$.
 Moreover, for any $p\ge 2$ we can show that positive-semidefiniteness of $g(\Re_{p}(\varphi),\,\varphi)$ for all $\varphi\in C^{\infty}(S_{0}^{p} M)$ and of $g(\Re_{p}(\varphi),\,\varphi)$ for all $\varphi \in C^{\infty}(S^{p} M)$ are equivalent (see \cite[p.~8]{17}).
In this case, we can reformulate Theorem~\ref{T-2.3} in the following form.

\begin{corollary}
Let $(M,g)$ be a complete noncompact Riemannian manifold with positive semi-definite sectional curvature and with the Lichnerowicz Laplacian
$\Delta_{L}:\,C^{\infty}(S^{p} M)\to C^{\infty}(S^{p} M)$ with $c>0$ acting on
$C^{\infty}$-sections of the bundle $S^{p} M$ of covariant symmetric $p$-tensor fields over $(M,g)$ for $p\ge 2$. Then the vector space $L^{q}({\rm Ker}\,\Delta_{L})$ is trivial for an arbitrary $q\in[1,+\infty)$.
\end{corollary}

The fact that ${\rm sec}\le0$ implies negative-semidefiniteness of $g(\Re_{p}(\varphi),\,\varphi)$ for all $p\ge2$
and $\varphi\in C^{\infty}(S_{0}^{p} M)$ was proved in \cite{DS,HMS}. Therefore, we can formulate a corollary from Theorem~\ref{T-2.4}.

\begin{corollary}
\hskip-0.4pt
Let $(M,g)$ be a complete noncompact Riemannian manifold~with negative semi-definite sectional curvature and $\Delta_{L}:C^{\infty}(S^{p} M)\to C^{\infty}(S^{p} M)$ be the Lichnerowicz Laplacian with $c<0$ acting on  $C^\infty$-sections of the bundle $S_0^p M$ of covariant symmetric traceless $p$-tensor fields over  $(M,g)$ for $p\ge2$.
Then the vector space $L^q({\rm Ker}\,\Delta_L)$ for any $q\in(0,+\infty)$ consists of parallel tensor fields.
In particular, if
 ${\rm vol}(M,g)= +\infty$, then $L^q({\rm Ker}\,\Delta_L)$ is trivial.
\end{corollary}

\textbf{3.2.} Second, we rewrite the Weitzenb\"{o}ck formula \eqref{GrindEQ-2-4} for the Lichnerowicz Laplacian $\Delta_{L}: C^{\infty}(S^{2} M)\to C^{\infty}(S^{2} M)$ in the following form (with $c\ne0$):
\begin{equation}\label{GrindEQ-3-1}
 \Delta_{L}\,\varphi =\bar{\Delta}\,\varphi +{c}\,\Re_{2} (\varphi).
\end{equation}
In this case,
the Weitzenb\"{o}ck curvature operator
\eqref{GrindEQ-2-5b} reduces to the form
\begin{eqnarray*}
 (\Re_{2}(\varphi))(X_1,X_2)&=& \sum\nolimits_{j}\big({\rm Ric}(e_j,X_1)\varphi(e_j,X_2) +{\rm Ric}(e_j,X_2)\varphi(e_j,X_1)\big) \\
 &-& 2\,\sum\nolimits_{j,k} R(e_j,X_1,e_k,X_2)\,\varphi(e_j,e_k) ,
\end{eqnarray*}
or, equivalently, \eqref{GrindEQ-2-5} has the following form (see \cite[p.~64]{2} and \cite{3,18}):
\begin{equation}\label{GrindEQ-3-2}
 (\Re_{2}(\varphi))_{ij} = R_{ik}\varphi_{j}^{k} +R_{jk}\varphi_{i}^{k}  -2R_{ikjl}\varphi^{kl}
\end{equation}
for the local components $\varphi_{ij}$ of an arbitrary $\varphi\in C^{\infty}(S^{2} M)$.
 Directly from \eqref{GrindEQ-3-1} and \eqref{GrindEQ-3-2} we obtain
${\rm trace}_{g}(\Delta_{L}\,\varphi )=\bar{\Delta}( {\rm trace}_{g}\,\varphi)$ for an arbitrary $\varphi \in C^{\infty}(S^{2} M)$. Therefore, the following statement holds (see also~\cite{1}).

\begin{proposition}
Let $\Delta_{L}:\,C^{\infty}(S^{2} M)\to C^{\infty}(S^{2} M)$ be the Lichnerowicz Laplacian acting on $C^{\infty}$-sections of the bundle ${S}^{2} M$ over a Riemannian manifold $({M,g})$, then
${\rm trace}_{g}(\Delta_{L}\,\varphi)=\bar{\Delta}( {\rm trace}_{g}\,\varphi )$.
\end{proposition}



It is well known that for any 2-tensor $\varphi$ the following inequality holds:
\[
 \|\varphi\|^2 \ge (1/n)({\rm trace}_g\,\varphi)^2
\]
at an arbitrary point $x\in M$. Therefore, if $\varphi\in L^2({\rm Ker}\,\Delta_{L})$, then ${\rm trace}_g\,\varphi\in L^2(M)$. On the other hand, we know that if $\varphi\in {\rm Ker}\,\Delta_{L}$, then $\overline{\Delta}({\rm trace}_g\,\varphi)=0$. At the same time, Yau proved in [38] that any harmonic function $f$ satisfying $f\in L^q(M)$ for some $q\in(1,+\infty)$ is constant on a complete manifold $(M,g)$. In particular, if ${\rm vol}(M,g)=+\infty$, then $f\equiv0$. In this case, from the above assumption we can conclude that the following holds.

\begin{proposition} Let $(M,g)$ be a complete Riemannian manifold and $\Delta_{L}: C^{\infty}(S^{2} M)\to C^{\infty}(S^{2} M)$ be the Lichnerowicz Laplacian. Then the trace of any smooth section of $L^2({\rm Ker}\,\Delta_{L})$ is a constant function. In particular, if the volume of $(M,g)$ is infinite, then $L^2({\rm Ker}\,\Delta_{L})$ consists of traceless symmetric 2-tensors.
\end{proposition}

 On the other hand, it is well known that there are no non-constant harmonic functions on a closed Riemannian manifold. Therefore, we can formulate the following corollary.

\begin{corollary}
Let $(M,g)$ be a closed Riemannian manifold with the Lichnerowicz Laplacian $\Delta_{L}: C^\infty(S^2 M)\to C^\infty(S^2 M)$ acting on $C^{\infty}$-sections of the bundle ${S}^{2} M$ over $(M,g)$. Then ${\rm trace}_{g}\,\varphi$ is a constant for an arbitrary bilinear form $\varphi\in {\rm Ker}\,\Delta_{L}$.
\end{corollary}

 In our case, \eqref{GrindEQ-2-6} can be rewritten in the following form  (with $c\ne0$):
\begin{equation}\label{GrindEQ-3-3}
 \frac{1}{2}\,\Delta_{B} (\|\varphi \|^{2}) = -g(\Delta_{L}\,\varphi,\,\varphi)
 +\| \nabla \varphi\|^{\; 2} + {c}\,g(\Re_{2} (\varphi),\,\varphi).
\end{equation}
Furthermore, for any point ${x}\in {M}$ there exists an orthonormal eigen-frame $\{e_{1},\ldots,e_{n}\}$ of $T_{x} M$ such that $\varphi_{x}(e_{i}, e_{j})=\mu_{i}\delta_{ij}$ for the Kronecker delta $\delta_{{ij}}$.
Then we have the formula (see \cite[p.~388]{3} and \cite[p.~436]{2})
\begin{equation}
 g(\Re_{2}(\varphi_{x}), \varphi_{x} ) = 2\sum\nolimits_{i<j}\sec(e_{i}\wedge e_{j})( \mu_{i} -\mu_{j})^{2},
\end{equation}
where $\sec({e}_{{i}}\wedge e_{j})={R}( {e}_{{i}}, e_{j}, e_{i}, e_{j})$ is the \textit{sectional curvature} $\sec(\sigma_{x})$ of $(M,g)$ in the direction of the tangent two-plane section
$\sigma_{{x}} ={\rm span}\{ e_{i}, e_{j}\}$ at ${x}\in{M}$.
Then we can rewrite \eqref{GrindEQ-3-3} in the following form:
\begin{equation*}
 \frac{1}{2}\,\Delta_{B} (\|\varphi \|^{2})
 =-g( \Delta_{L}\,\varphi, \varphi) +\|\nabla \varphi \|^{2}
 + 2\,c\sum\nolimits_{i<j}\sec(e_{i}\wedge e_{j} )(\mu_{i} -\mu_{j})^{2}  .
\end{equation*}
In particular, if $\varphi$ is a covariant $\Delta_{L}$-harmonic 2-tensor, then we have
\begin{equation}\label{GrindEQ-3-6}
 \frac{1}{2}\,\Delta_{B}(\|\varphi \|^{2}) = \|\nabla \varphi \|^{2}
 +2\,c\sum\nolimits_{i<j}\sec(e_{i}\wedge e_{j}) (\mu_{i} -\mu_{j} )^{2}  .
\end{equation}
From \eqref{GrindEQ-3-6} we conclude that $\|\varphi\|^{2}$ is a subharmonic function if $c>0$ and the sectional curvature of $(M,g)$ is non-negative. Therefore, proceeding from the above formula and using the Hopf maximum principle (see \cite[p.~26]{6} and \cite{11}), we can conclude the following: if the sectional curvature $\sec(\sigma_{x})$ of $(M,g)$ is non-negative at any point ${x}$ of a connected open domain ${U}\subset {M}$ and $\sec(\sigma_{x})$ is strictly positive (in all 2-dimensional directions $\sigma_{{x}} $) at some point $x$ of ${U}$, then $\|\varphi\,\|^{2}$ is a constant $C$ and $\nabla\varphi =0$ in ${U}$. If $C>0$, then $\varphi$ is nowhere zero.
Now, at a point $x$ of ${U}$, where the sectional curvature $\sec(\sigma_{x})$ is positive, the left hand side of \eqref{GrindEQ-4-6} is zero, while the right hand side is nonpositive. This contradiction shows $\mu_{1} =\ldots=\mu_{n} =\mu$ and hence $\varphi =\mu \cdot g$ for some constant $\mu$ everywhere in ${U}$. On the other hand, the fact that $\nabla\,\varphi =0$ means that $\varphi$ is invariant under parallel translation. In this case, if the holonomy of $(M,g)$ is irreducible, then the tensor $\varphi$ has a one eigenvalue, i.e., $\varphi =\mu \cdot {g}$ for some constant $\mu$ at each point of ${U}$. As a result, we have the following local theorem.

\begin{theorem}\label{T-3.1}
Let ${U}$ be a connected open domain of a Riemannian manifold $(M,g)$ with nonnegative sectional curvature at any point of $U$ and $\Delta_{L}: C^{\infty}(S^{2} M)\to C^{\infty}(S^{2} M)$ the Lichnerowicz Laplacian with $c>0$ acting on $C^{\infty}$-sections of the bundle ${S}^{2} M$ over $({M,g})$. If $\varphi \in {\rm Ker}\,\Delta_{L}$ at any point of ${U}$ and the scalar function $\|\varphi\,\|^{2}$ has a local maximum at some point of ${U}$, then $\|\varphi\,\|^{2}$ is a constant function and $\varphi$ is invariant under parallel translation in ${U}$. Moreover, if $\sec(\sigma_{x}) >0$ in all directions $\sigma_{{x}}$ at some point $x\in U$ or the holonomy of $(M,g)$ is irreducible, then $\varphi$ is a constant multiple of $g$ at each point of ${U}$.
\end{theorem}

Based on \eqref{GrindEQ-3-3} and using our Theorem~\ref{T-2.3} and the Greene-Wu theorem on subharmonic functions on a complete noncompact Riemannian manifold $(M,g)$ with nonnegative sectional curvature (see \cite{14}), we conclude that if $\varphi \in {\rm Ker}\,\Delta_{L}$ at any point of $(M,g)$ and $\int_{M}\,\|\varphi\|^{\,2}\,{\rm d}\,{\rm vol}_g <+\infty$, then $\varphi \equiv 0$. Then we can formulate the following.

\begin{corollary}
Let $(M,g)$ be a complete noncompact Riemannian manifold with nonnegative sectional curvature and $\Delta_{L}: C^{\infty}(S^{2} M)\to C^{\infty}(S^{2} M)$ the Lichnerowicz Laplacian with $c>0$ acting on $C^{\infty}$-sections of the bundle ${S}^{2} M$ over $(M,g)$. Thus the vector space $L^{2}({\rm Ker}\,\Delta_{L})$ is trivial.
\end{corollary}

 Consider now the case $n=3$. We have the following equality:
\[
 \sec(\sigma_{{x}})= (1/2) s - {\rm Ric}(X_{x},X_{x}),
\]
where $\sec(\sigma_{x})$ is the sectional curvature in the direction of the plane $\sigma_{x}\subset T_{x} M$ for an arbitrary point $x\in M^3$, ${X}$ is a unit vector orthogonal to $\sigma_{x}$, and $s$ is the scalar curvature of $(M, g)$
 (see \cite[Lemma~2.1]{19}). Therefore, if $n=3$ and ${\rm Ric}\le (1/2) s\,g$ at each point $x\in M$, then the inequality $\sec(\sigma_{x})\ge 0$ holds at each point $x\in M$. In this case, if $c>0$ then from \eqref{GrindEQ-3-6} we conclude that $\|\varphi \|^{\, 2}$ is a subharmonic function, and using the Greene-Wu theorem on subharmonic functions, we conclude that $\varphi \equiv 0$. Thus, we can formulate the following.

\begin{corollary}
Let $(M,g)$ be a three-dimensional complete noncompact Riemannian manifold and
$\Delta_{L}: C^{\infty}(S^{2} M)\to C^{\infty}(S^{2} M)$ the Lichnerowicz Laplacian acting on $C^{\infty}$-sections of the bundle ${\rm S}^{2} M$ over $(M,g)$.
If the Ricci curvature ${\rm Ric}$ and the scalar curvature $s$ of $(M, g)$ satisfy the inequality ${\rm Ric}\le(1/2)s\,g$, then the vector space $L^{2}( {\rm Ker}\,\Delta_{L})$ is trivial.
\end{corollary}

\textbf{3.3.} Consider a closed Riemannian manifold $(M,g)$ with nonnegative sectional curvature. Then, based on \eqref{GrindEQ-3-6} and the Bochner maximum principle (see \cite[p.~30]{6}), we can conclude that the kernel of the Lichnerowicz Laplacian $\Delta_{L}:\,C^{\infty}(S^{2} M)\to C^{\infty}(S^{2} M)$ with $c>0$ consists of parallel symmetric 2-tensor tensor fields on $(M, g)$, i.e., from the condition $\varphi\in{\rm Ker}\,\Delta_{L}$ we obtain
$\nabla\,\varphi =0$. This equation means that
\[
 \varphi_{im} R_{jkl}^{m} + \varphi_{jm} R_{ikl}^{m} =0.
\]
The last algebraic equalities imply that $\varphi_{i}^{k} R_{kj} -R_{ikjl}\,\varphi^{kl} =0$.
Combining this result with the Bianchi identities, we obtain
\begin{eqnarray*}
 (\Re_{2}(\varphi))_{ij} = R_{ik}\varphi_{j}^{k} +R_{jk}\varphi_{i}^{k}
 -2R_{ikjl}\varphi^{kl} \\
 =( R_{ik}\varphi_{j}^{k} -R_{ikjl}\varphi^{kl}) +( R_{jk}\varphi_{i}^{k} -R_{jkil}\varphi^{kl} ) =0,
\end{eqnarray*}
where $R_{ikjl}\varphi^{kl} =(-R_{ijlk} -R_{ilkj}) \varphi^{kl} =-R_{ilkj}\,\varphi^{kl}
=R_{jkil}\varphi^{kl}$.
In this case, from \eqref{GrindEQ-3-1}
we obtain $\bar{\Delta}\,\varphi_{ij} =0$, i.e., $\varphi_{ij}$ is a harmonic function for any $i,j=1,\ldots, n$.
It is well known that a harmonic function on a closed Riemannian manifold is constant (see \cite[p.~30]{6}). Therefore, we can formulate the following theorem.

\begin{theorem}
Let $(M, g)$ be a closed Riemannian manifold with nonnegative sectional curvature. Then the kernel of the Sampson Laplacian $\Delta_{S}: C^{\infty}(S^{2} M)$ $\to C^{\infty}(S^{2} M)$ with $c>0$ consists of constant symmetric 2-tensors.
\end{theorem}

Note that a \textit{Riemannian symmetric space of compact type} provides an example of a closed Riemannian manifold with non-negative sectional curvature and positive-definite Ricci tensor (see \cite[p.~256]{15}).
Thus the following corollary is true.

\begin{corollary}
Let $(M,g)$ be a Riemannian symmetric space $(M,g)$ of compact type.
Then the kernel of the Lichnerowicz Laplacian $\Delta_{L}: C^{\infty}(S^{2} M)\to C^{\infty}(S^{2} M)$ with $c>0$ consists of constant covariant symmetric 2-tensors.
\end{corollary}

\begin{remark}\rm
A simple example of Riemannian symmetric spaces of compact type is the ${n}$-dimensional round sphere $(S^{n}, {g}_{0}$) with standard metric ${g}_{0}$.
Then an arbitrary $\Delta_{L}$-harmonic tensor on $(S^{n}, {g}_{0}$) has the form $\varphi=\mu \cdot g_{0}$ for some real constant $\mu$.
\end{remark}

In conclusion, we recall the definition of a \textit{TT-tensor} (Transverse Traceless tensor), that is a divergence free and traceless covariant symmetric 2-tensor field. Such tensors are of fundamental importance in \textit{stability analysis} in General Relativity (see, for example, \cite{22,24,23}) and in Riemannian geometry (see \cite{25,1}). In particular, Page and Pope have proved in \cite{24} the following theorem on the kernel of the Lichnerowicz Laplacian acting on \textit{TT}-tensors.

\begin{theorem}
Let $(M,\,g)$ be a connected Riemannian manifold and $\Delta_{L}$ the Lichnerowicz Laplacian with $c=1$ acting on $C^{\infty}$-sections of the bundle $S^{2} M$ over $(M,\,g)$. If the holonomy of $(M,\,g)$ is reducible, then there exists a \textit{TT}-tensor $\varphi \in C^{\infty}(S^{2} M)$ such that $\varphi \in {\rm Ker}\,\Delta_{L}$.
\end{theorem}

 Side by side, one can prove the following corollary of our Theorem~\ref{T-3.1} for $\Delta_{L}$-harmonic \textit{TT}-tensors.

\begin{corollary}
Let $(M,g)$ be a closed Riemannian manifold with positive sectional curvature and $\Delta_{L}:\,C^{\infty}(S^{2} M)\to C^{\infty}(S^{2} M)$ the Lichnerowicz Laplacian with $c=1$ restricted to $TT$-tensors on $(M,g)$. Then the vector space $L^{2}( {\rm Ker}\,\Delta_{L})$ is trivial.
\end{corollary}

It is known that in dimension three a metric $g$ has positive sectional curvature if and only if its Ricci curvature ${\rm Ric}$ and scalar curvature $s$ satisfy the inequality ${\rm Ric}<(1/2)s\,g$ (see \cite[p.~277]{26}).
Therefore, we can formulate the following.

\begin{corollary}
Let $(M,g)$ be a closed Riemannian manifold and $\Delta_{L}{:}\,C^{\infty}(S^{2} M)$ $\to C^{\infty}(S^{2} M)$ the Lichnerowicz Laplacian with $c=1$ restricted to \textit{TT}-tensors on $(M, g)$. If~the Ricci curvature ${\rm Ric}$ and the scalar curvature $s$ of $(M,g)$ satisfy the inequa\-lity ${\rm Ric} < (1/2)\,s\,g$, then the vector space $L^{2}({\rm Ker}\,\Delta_{L})$ is trivial.
\end{corollary}

\section{\large Applications to the theories of infinitesimal Einstein deformations
and the stability of Einstein manifolds}

For the case $\Delta_{L}:\,C^{\infty}(S^{2} M)\to C^{\infty}(S^{2} M)$, the Lichnerowicz Laplacian $\Delta_{L}$ with $c=1$ is of fundamental importance in the stability analysis in General Relativity (see, for instance, \cite{17,22,23}) and appears in many problems of Riemannian geometry. For example, the Lichnerowicz Laplacian acting on symmetric 2-tensor fields can be seen as infinitesimal deformations of metric $g$, and describes the change of the Ricci tensor
in terms of these infinitesimal deformations (see, for example, \cite{3} and \cite[Chapter~12]{2}). Furthermore, the Lichnerowicz Laplacian is a fundamental operator; when acting on covariant symmetric 2-tensor fields in context of \textit{Ricci flow}, it seems to be more natural than the rough Laplacian $\bar{\Delta}$. Examples of this naturality are the appearance of $\Delta_{L}$ in the linearized Ricci flow equation (see, for example, the evolution formula of the Ricci tensor under the Ricci flow in \cite[p.~112]{27}).
 In this section, we complete these results.

\smallskip

\textbf{4.1.} Recall that an Einstein manifold is an $n$-dimensional Riemannian manifold $(M, g)$, for which the Ricci tensor satisfies ${\rm Ric}=\kappa\,g$ for some real number~$\kappa$. Taking trace of this, one can prove that $\kappa ={s / n}$ for the scalar curvature $s$ of $(M,\,g)$. We shall consider Einstein manifolds in this section.

Notice that the Riemannian curvature tensor of $(M, g)$ defines a symme\-tric algebraic operator ${\mathop{R}\limits^{\circ}}:\,S^{2}(T_{x} M)\to S^{2}(T_{x} M)$ on the vector space $S^{2}(T_{x} M)$ of symmetric bilinear forms over tangent space $T_{x} M$ at an arbitrary point $x\in M$.
The~operator ${\mathop{R}\limits^{\circ}}$ is called the \textit{curvature operator of the second kind} of $(M, g)$.

\begin{remark}\rm
The definition, properties and applications of ${\mathop{R}\limits^{\circ}}$ can be found in monographs \cite{2,7} and in articles from the following list: \cite{29,25,30,31,21,28,32}.
\end{remark}

Then we call the differential operator
\[
 \Delta_{E} =\bar{\Delta}-2\,{\mathop{R}\limits^{\circ}},
\]
acting on $C^{\infty}$-sections of the bundle $S^{2} M$ over an Einstein manifold $(M,\,g)$ the \textit{Einstein operator}.
This is a self-adjoint elliptic operator mapping from the vector space of \textit{TT}-tensors to itself (see also \cite{33}).
If a \textit{TT}-tensor $\varphi$ belongs to ${\rm Ker}\,\Delta_{E}$ then it can be seen as an \textit{infinitesimal Einstein deformation through} $g$ (see \cite{3} and \cite[pp.~346--348]{2}). Recall that a deformation of Einstein structures through $g$ means a smooth curve $g(t)$ of Riemannian metrics, where $t$ belongs to some open interval containing $0$ with $g(0)=g$ and such that for each $t$ there exists a real number $\kappa(t)$ with the property ${\rm Ric}_{g(t)}=\kappa(t)\cdot g(t)$.

The Einstein operator is closely related to the Lichnerowicz Laplacian $\Delta_{L}$ with $c=1$ .
In fact, on Einstein manifolds, we have the relation
\begin{equation}\label{GrindEQ-4-1}
 \Delta_{L}=\Delta_{E} +2\,({s}/{n})\,{\rm Id}.
\end{equation}
Therefore, if $\varphi \in C^{\infty}(S^{2} M)\cap {\rm Ker}\,\Delta_{L}$ then $\Delta_{E}\,\varphi =-2({s}/{n})\,\varphi$, i.e., $\varphi$ is an eigentensor of $\Delta_{E}$ with the eigenvalue $-2\,{s / n}$.
The converse is also true. From \eqref{GrindEQ-4-1} we can deduce that the Einstein operator $\Delta_{E}$ is positive (resp., negative) for all \textit{TT}-tensors belonging to ${\rm Ker}\,\Delta_{L}$,
if $(M,g)$ is an Einstein manifold with negative (resp., positive) scalar curvature. In particular, if $(M,\,g)$ is a \textit{Ricci-flat} Riemannian manifold (see \cite{34}), then $\Delta_{L}=\Delta_{E}$. In this case, an arbitrary \textit{TT}-tensor $\varphi$ is an infinitesimal Einstein deformation of the metric $g$ if $\varphi$ belongs to ${\rm Ker}\,\Delta_{L}$. Therefore, we can formulate the following.

\begin{proposition}\label{P-4.1}
Let $(M,g)$ be an Einstein manifold, then ${\rm Ker}\,\Delta_{L}$ of the Lichnerowicz Laplacian $\Delta_{L}: C^{\infty}(S^{2} M)\to C^{\infty}(S^{2} M)$ consists of eigentensors of the Einstein operator $\Delta_{E} =\bar{\Delta}-2\,{\mathop{R}\limits^{\circ}}$ with eigenvalues equal to $-2\,s /n$. The converse is also true. Furthermore, the Einstein operator $\Delta_{{E}}$ is positive (resp., negative) on \textit{TT}-tensors belonging to ${\rm Ker}\,\Delta_{L}$, if the scalar curvature is negative (resp., positive). In particular, if $(M,g)$ is a Ricci-flat Riemannian manifold, then an arbitrary \textit{TT}-tensor $\varphi$ belongs to ${\rm Ker}\,\Delta_{L}$ if and only if it is an infinitesimal Einstein deformation of the metric $g$.
\end{proposition}

Recall that $(M,g)$ is called \textit{unstable}, if the Einstein operator admits negative eigenvalues on \textit{TT}-tensors (see \cite{33}).
From our Proposition \ref{P-4.1} we obtain that $(M,g)$ is unstable with respect to a \textit{TT}-tensor $\varphi \in{\rm Ker}\,\Delta_{L} $, if $(M,g)$ is an Einstein manifold with positive scalar curvature.

The following theorem on infinitesimal Einstein deformations from \cite[p.~355]{2} is well known. Let ${g}$ be an Einstein metric on $M$ and by $a_{0}$ -- the largest eigenvalue of the zero order operator ${\mathop{R}\limits^{\circ}}$ on the bundle of trace-free symmetric 2-tensor fields, i.e.,
\[
 {a}_{0}={\rm sup}\{\,g({\mathop{R}\limits^{\circ}} h, h) / \|h\|^{2}:\ h\in C^{\infty}(S_{0}^{2} M)\,\}.
\]
If $a_{0} <\max \{- s / n; s /(2n)\}$, then $g$ has no infinitesimal Einstein deformations.
On the other hand, the following our theorem completes this theorem.

\begin{theorem}
Let $(M,g)$ be a closed Einstein manifold with nonzero scalar curvature $s$
and ${K}_{\rm min}$ -- the minimum of its sectional curvature.
If ${K}_{\rm min}\ge {s}/{n^{2}}$, then $(M, g)$ is not an unstable manifold and does not admit infinitesimal Einstein deformations.
\end{theorem}

\noindent\textbf{Proof}.
Let $(M,g)$ be an Einstein manifold with nonzero scalar curvature $s$ and let $\varphi$ be a \textit{TT}-tensor on $(M,g)$. Then \eqref{GrindEQ-3-3} can be rewritten in the form
\begin{eqnarray}\label{GrindEQ-4-3}
\nonumber
 \frac{1}{2}\,\Delta_{B}(\|\,\varphi\,\|^{2}) =-{g}(\Delta_{{E}}\,\varphi,\varphi)
 - 2\,\frac{s}{n}\,\|\varphi\|^{2}  + \| \nabla  \varphi \|^{2}\\
 +\, 2 \sum\nolimits_{i<j}\sec(e_{i}\wedge e_{j})(\mu_{i} -\mu_{j} )^{2} ,
\end{eqnarray}
where $c=1$.
If ${\rm trace}_{g}\,\varphi =\mu_{1} + \ldots +\mu_{n} = 0$,
then
the following equality holds:
\[
 \|\varphi\|^{2} =\mu_{1}^{2}
 +\ldots+\mu_{n}^{2} =\frac{1}{n}\sum\nolimits_{i<j}(\mu_{i} -\mu_{j})^{2}  .
\]
In this case, from \eqref{GrindEQ-4-3} one can obtain the inequality
\begin{equation}\label{GrindEQ-4-4}
 \frac{1}{2}\,\Delta_{B}(\|\varphi \|^{2}) \ge -{g}(\Delta_{{E}}\,\varphi, \varphi)
 +\left\|\nabla \varphi \right\|^{2} +2\,\big({K}_{\rm min} -\frac{s}{n^{2}}\big) \sum\nolimits_{i<j}(\mu_{i} -\mu_{j})^{2},
\end{equation}
where we denoted by ${K}_{\rm min}$ the minimum of the sectional curvature of $(M,g)$, i.e., $\sec(\sigma_{x})\ge {K}_{\rm min}$ in all 2-dimensional directions $\sigma_{{x}}$ at each point $x\in M$.

 First, let $(M,g)$ be unstable with respect to $\varphi $, then
\[
 {g}(\Delta_{{E}}\,\varphi, \varphi ) = -\lambda^{2}(\varphi)\|\varphi\,\|^{2}
\]
 for some $\lambda (\varphi)\ne 0$. In this case, the inequality \eqref{GrindEQ-4-4} takes the form
\begin{equation}\label{GrindEQ-4-5}
 \frac{1}{2}\,\Delta_{B}(\,\|\varphi\|^{2}) \ge \lambda^{2}(\varphi)
 \|\varphi \|^{2} + \|\nabla \varphi \|^{\, 2}
 + 2\,\big({K}_{\rm min} -\frac{s}{n^{2}}\big)\sum\nolimits_{i<j}(\mu_{i} -\mu_{j} )^{2} .
\end{equation}
If ${K}_{\rm min}\ge {s}/{n^{2}} $, then from \eqref{GrindEQ-4-5} we conclude that $\|\varphi\|^{\, 2}$ is a subharmonic function, i.e., $\Delta_{B}(\|\varphi\|^{\, 2})\ge 0$. Furthermore, if $(M,g)$ is a closed manifold, then using the Bochner maximum principle (see \cite[p.~30]{6}), we conclude that $\|\varphi\|^{2}$ is constant.
In this case, from \eqref{GrindEQ-4-5} we obtain that $\varphi \equiv 0$.

 Second, let $(M,g)$ be a stable manifold, then \eqref{GrindEQ-4-4} can be rewritten in the form
\begin{equation}\label{GrindEQ-4-6}
 \frac{1}{2}\,\Delta_{B}\,(\|\,\varphi\,\|^{\, 2}) \ge \|\,\nabla\,\varphi\,\|^{\, 2}
 + 2\,\big({K}_{\rm min} -\frac{s}{n^{2}}\big)\sum\nolimits_{i<j}(\mu_{i} -\mu_{j})^{2}  .
\end{equation}
If ${K}_{\rm min}\ge {s}/{n^{2}}$, then from \eqref{GrindEQ-4-6} we obtain $\Delta_{B}(\|\varphi\|^{\,2})\ge 0$, i.e., $\Delta_{B}(\|\varphi \|^{\, 2})$ is a subharmonic function. Then proceeding from \eqref{GrindEQ-4-6} and using the Bochner maximum principle (see \cite[p.~30]{6}), we conclude that $\|\varphi\|={\rm const}$ and hence $\nabla\,\varphi =0$. In this case, by the Ricci identities, we have $\varphi_{ik} R^{k}_{jlm} +\varphi_{kj} R^{k}_{ilm} =0$.
Then ${\mathop{R}\limits^{\circ\,}} (\varphi)=-({s}/{n})\varphi$. In this case, the equation $\Delta_{E}\,\varphi =\bar{\Delta}\,\varphi -2\,{\mathop{R}\limits^{\circ}}\,(\varphi)=0$ can be rewritten in the form $\bar{\Delta}\,\varphi =- 2({s}/{n})\varphi$. This implies
\[
 -2\,\frac{s}{n}\int_{M}\|\varphi \|^{2}\,{\rm d}\,{\rm Vol}_{g} =\int_{M}g( \bar{\Delta}\varphi, \varphi )\,{\rm d}\,{\rm vol}_{g} = \int_{M}\|\nabla \varphi \|^{2}\,{\rm d}\,{\rm vol}_{g} =0.
\]
Hence, $\varphi \equiv 0$. By this, $\varphi$ is a trivial infinitesimal Einstein deformation.
\hfill$\Box$

\end{document}